\documentclass[11pt]{article}

\usepackage{amsmath,amssymb}
\usepackage{a4wide}

\begin{document}

\title{Consecutive, Reversed, Mirror, and Symmetric\\
Smarandache Sequences of Triangular Numbers\footnote{Research
Report CM04/D-05, June 2004. Submitted to the Smarandache Notions Journal.}}

\author{Delfim F. M. Torres \\
        \texttt{delfim@mat.ua.pt} \\
        Department of Mathematics \\
        University of Aveiro, Portugal
        \and
        Viorica Teca\footnote{MSc student
        of Informatics at the University of Craiova, Romania.
        Student at the University of Aveiro under the
        Socrates/Erasmus European programme, 2004.} \\
        \texttt{viorica\_teca@yahoo.com} \\
        Faculty of Mathematics-Informatics \\
        University of Craiova, Romania
        }

\date{07/June/2004}

\maketitle


\begin{abstract}
\noindent We use the \textsf{Maple} system to check the  investigations of S.~S.~Gupta
regarding the Smarandache consecutive and the reversed Smarandache sequences
of triangular numbers [Smarandache Notions Journal, Vol. 14, 2004, pp.~366--368].
Furthermore, we extend previous investigations to the mirror
and symmetric Smarandache sequences of triangular numbers.
\end{abstract}


\bigskip


\bigskip


The $n$th triangular number $t_n$,  $n \in \mathbb{N}$,
is defined by $t_n = \sum_{i=1}^{n} i = n(n+1)/2$.
These numbers were first studied by the Pythagoreans.

The first $k$ terms of the triangular sequence $\{t_n\}_{n=1}^{\infty}$
are easily obtained in \textsf{Maple}:
\begin{verbatim}
  > t:=n->n*(n+1)/2:
  > first := k -> seq(t(n),n=1..k):
\end{verbatim}
For example:
\begin{verbatim}
  > first(20);
\end{verbatim}
\vspace{-0.3cm}
\begin{equation*}
1, 3, 6, 10, 15, 21, 28, 36, 45, 55, 66, 78, 91, 105, 120, 136, 153,
171, 190, 210
\end{equation*}
In this short note we are interested in studying Smarandache
sequences of triangular numbers with the help of the
\textsf{Maple} system.

To define the Smarandache sequences, it is convenient to introduce
first the concatenation operation. Given two positive integer
numbers $n$ and $m$, the concatenation operation \texttt{conc} is defined
in \textsf{Maple} by the following function:
\begin{verbatim}
  > conc := (n,m) -> n*10^length(m)+m:
\end{verbatim}
For example,
\begin{verbatim}
  > conc(12,345);
\end{verbatim}
\vspace{-0.3cm}
\begin{equation*}
12345
\end{equation*}
Given a positive integer sequence $\{u_n\}_{n=1}^{\infty}$,
we define the corresponding \emph{Smarandache Consecutive Sequence}
$\{scs_n\}_{n=1}^{\infty}$ recursively:
\begin{equation*}
\begin{split}
scs_1 &= u_1 \, ,\\
scs_n &= conc(scs_{n-1},u_n) \, .
\end{split}
\end{equation*}
In \textsf{Maple} we define:
\begin{verbatim}
  > scs_n := (u,n) -> if n=1 then u(1) else conc(scs_n(u,n-1),u(n)) fi:
  > scs := (u,n) -> seq(scs_n(u,i),i=1..n):
\end{verbatim}
The standard Smarandache consecutive sequence,
introduced by the the Romanian mathematician Florentin Smarandache,
is obtained when one chooses $u_n = n$, $\forall n \in \mathbb{N}$.
The first 10 terms are:
\begin{verbatim}
  > scs(n->n,10);
\end{verbatim}
\vspace{-0.3cm}
\begin{equation*}
  1, 12, 123, 1234, 12345, 123456, 1234567, 12345678, 123456789,
  12345678910
\end{equation*}
Another example of a Smarandache consecutive sequence
is the Smarandache consecutive sequence of triangular numbers.
With our \textsf{Maple} definitions, the first 10 terms of such sequence are
obtained with the following command:
\begin{verbatim}
  > scs(t,10);
\end{verbatim}
\vspace{-0.3cm}
\begin{gather*}
1, 13, 136, 13610, 1361015, 136101521, 13610152128, \\
1361015212836, 136101521283645, 13610152128364555
\end{gather*}
Sometimes, it is preferred to display Smarandache sequences
in ``triangular form''.
\begin{verbatim}
  > show := L -> map(i->print(i),L):
  > show([scs(t,10)]):
\end{verbatim}
\vspace{-0.3cm}
\begin{gather*}
                                  1\\
                                  13\\
                                 136\\
                                13610\\
                               1361015\\
                              136101521\\
                             13610152128\\
                            1361015212836\\
                           136101521283645\\
                          13610152128364555
\end{gather*}

The \emph{Reversed Smarandache Sequence} ($rss$) associated
with a given sequence $\{u_n\}_{n=1}^{\infty}$, is defined recursively by
\begin{equation*}
\begin{split}
rss_1 &= u_1 \, ,\\
rss_n &= conc(u_n,rss_{n-1}) \, .
\end{split}
\end{equation*}
In \textsf{Maple} we propose the following definitions:
\begin{verbatim}
  > rss_n := (u,n) -> if n=1 then u(1) else conc(u(n),rss_n(u,n-1)) fi:
  > rss := (u,n) -> seq(rss_n(u,i),i=1..n):
\end{verbatim}
The first terms of the reversed Smarandache sequence
of triangular numbers are now easily obtained:
\begin{verbatim}
  > rss(t,10);
\end{verbatim}
\vspace{-0.3cm}
\begin{gather*}
  1, 31, 631, 10631, 1510631, 211510631, 28211510631, \\
  3628211510631, 453628211510631, 55453628211510631
\end{gather*}

We define the \emph{Smarandache Mirror Sequence} ($sms$) as follows:
\begin{equation*}
\begin{split}
sms_1 &= u_1 \, , \\
sms_n &= conc(conc(u_n,sms_{n-1}),u_n)
\end{split}
\end{equation*}
\begin{verbatim}
  > sms_n := (u,n) -> if n=1 then
  >                     u(1)
  >                   else
  >                     conc(conc(u(n),sms_n(u,n-1)),u(n))
  >                   fi:
  > sms := (u,n) -> seq(sms_n(u,i),i=1..n):
\end{verbatim}
The first 10 terms of the Smarandache mirror sequence introduced by Smarandache are:
\begin{verbatim}
  > sms(n->n,10);
\end{verbatim}
\vspace{-0.3cm}
\begin{gather*}
  1, 212, 32123, 4321234, 543212345, 65432123456, 7654321234567,\\
  876543212345678, 98765432123456789, 109876543212345678910
\end{gather*}
We are interested in the Smarandache mirror sequence of triangular numbers.
The first 10 terms are:
\begin{verbatim}
  > sms(t,10);
\end{verbatim}
\vspace{-0.3cm}
\begin{gather*}
  1, 313, 63136, 106313610, 1510631361015, 21151063136101521,\\
  282115106313610152128, 3628211510631361015212836, \\
  45362821151063136101521283645, 554536282115106313610152128364555
\end{gather*}
Finally, we define the Smarandache Symmetric Sequence ($sss$).
For that we introduce the function ``But Last Digit'' (\texttt{bld}):
\begin{verbatim}
  > bld := n -> iquo(n,10):
  > bld(123);
\end{verbatim}
\vspace{-0.3cm}
\begin{equation*}
12
\end{equation*}
If the integer number is a one-digit number,
then function \texttt{bld} returns zero:
\begin{verbatim}
  > bld(3);
\end{verbatim}
\vspace{-0.3cm}
\begin{equation*}
0
\end{equation*}
This is important: with our \texttt{conc} function,
the concatenation of zero with a positive integer $n$ gives $n$.
\begin{verbatim}
  > conc(bld(1),3);
\end{verbatim}
\vspace{-0.3cm}
\begin{equation*}
3
\end{equation*}
The \emph{Smarandache Symmetric Sequence} ($sss$) is now easily defined,
appealing to the Smarandache consecutive, and reversed Smarandache sequences:
\begin{equation*}
\begin{split}
sss_{2n-1} &= conc(bld(scs_{2n-1}),rss_{2n-1}) \, , \\
sss_{2n} &= conc(scs_{2n},rss_{2n}) \, ,
\end{split}
\end{equation*}
$n \in \mathbb{N}$. In \textsf{Maple}, we give the following definitions:
\begin{verbatim}
  > sss_n := (u,n) -> if type(n,odd) then
  >                     conc(bld(scs_n(u,(n+1)/2)),rss_n(u,(n+1)/2))
  >                   else
  >                     conc(scs_n(u,n/2),rss_n(u,n/2))
  >                   fi:
  > sss := (u,n) -> seq(sss_n(u,i),i=1..n):
\end{verbatim}
The first terms of Smarandache's symmetric sequence are
\begin{verbatim}
  > sss(n->n,10);
\end{verbatim}
\vspace{-0.3cm}
\begin{equation*}
1, 11, 121, 1221, 12321, 123321, 1234321, 12344321, 123454321, 1234554321
\end{equation*}
while the first 10 terms of the Smarandache symmetric sequence
of triangular numbers are
\begin{verbatim}
  > sss(t,10);
\end{verbatim}
\vspace{-0.3cm}
\begin{equation*}
1, 11, 131, 1331, 13631, 136631, 136110631, 1361010631,
1361011510631, 13610151510631
\end{equation*}

One interesting question is to find prime numbers in the above
defined Smarandache sequences of triangular numbers.
We will restrict our search to the first 1000 terms of each sequence.
All computations were done with \textsf{Maple} 9 running on a 2.00Ghz
\textsf{Pentium} 4 with 256Mb RAM.

We begin by collecting four lists with the first 1000 terms of
the consecutive, reversed, mirror, and symmetric Smarandache
sequences of triangular numbers:
\begin{verbatim}
  > st:=time(): Lscs1000:=[scs(t,1000)]: printf("%a seconds",round(time()-st));
    20 seconds
\end{verbatim}
\begin{verbatim}
  > st:=time(): Lrss1000:=[rss(t,1000)]: printf("%a seconds",round(time()-st));
    75 seconds
\end{verbatim}
\begin{verbatim}
  > st:=time(): Lsms1000:=[sms(t,1000)]: printf("%a seconds",round(time()-st));
    212 seconds
\end{verbatim}
\begin{verbatim}
  > st:=time(): Lsss1000:=[sss(t,1000)]: printf("%a seconds",round(time()-st));
    26 seconds
\end{verbatim}
We note that $scs_{1000}$ and $rss_{1000}$ are positive integer
numbers with 5354 digits;
\begin{verbatim}
  > length(Lscs1000[1000]), length(Lrss1000[1000]);
\end{verbatim}
\vspace{-0.3cm}
\begin{equation*}
5354, 5354
\end{equation*}
while $sms_{1000}$ and $sss_{1000}$ have, respectively,
10707 and 4708 digits.
\begin{verbatim}
  > length(Lsms1000[1000]), length(Lsss1000[1000]);
\end{verbatim}
\vspace{-0.3cm}
\begin{equation*}
10707, 4708
\end{equation*}
There exist two primes (13 and 136101521)
among the first 1000 terms of the
Smarandache consecutive sequence of triangular numbers;
\begin{verbatim}
  > st := time():
  > select(isprime,Lscs1000);
  > printf("%a minutes",round((time()-st)/60));
\end{verbatim}
\vspace{-0.3cm}
\begin{equation*}
[13, 136101521]
\end{equation*}
\vspace{-0.5cm}
\begin{verbatim}
    9 minutes
\end{verbatim}
six primes among the first 1000 terms of the reversed Smarandache
sequence of triangular numbers;
\begin{verbatim}
  > st := time():
  > select(isprime,Lrss1000);
  > printf("%a minutes",round((time()-st)/60));
\end{verbatim}
\vspace{-0.3cm}
\begin{equation*}
[31, 631, 10631, 55453628211510631, 786655453628211510631, 10591786655453628211510631]
\end{equation*}
\vspace{-0.5cm}
\begin{verbatim}
    31 minutes
\end{verbatim}
only one prime (313) among the first 600 terms of the Smarandache
mirror sequence of triangular numbers;\footnote{Our computer
runs low in memory when one tries the first 1000 terms of the
Smarandache mirror sequence of triangular numbers.
For this reason, we have considered here only the first 600 terms of the sequence.}
\begin{verbatim}
  > length(Lsms1000[600]); # sms_{600} is a number with 5907 digits
\end{verbatim}
\vspace{-0.3cm}
\begin{equation*}
5907
\end{equation*}
\begin{verbatim}
  > st := time():
  > select(isprime,Lsms1000[1..600]);
  > printf("%a minutes",round((time()-st)/60));
\end{verbatim}
\vspace{-0.3cm}
\begin{equation*}
[313]
\end{equation*}
\vspace{-0.5cm}
\begin{verbatim}
    3 minutes
\end{verbatim}
and five primes among the first 1000 terms of the
Smarandache symmetric sequence of triangular numbers
(the fifth prime is an integer with 336 digits).
\begin{verbatim}
  > st := time():
  > select(isprime,Lsss1000);
  > printf("%a minutes",round((time()-st)/60));
\end{verbatim}
\vspace{-0.3cm}
\begin{gather*}
  [11, 131, 136110631, 1361015212836455566789110512012010591786655453628211510631,\\
1361015212836455566789110512013615317119021023125327630032535137840643546549652856159\\
5630666703741780820861903946990103510811128117612251275132613781431148515401596165316\\
5315961540148514311378132612751225117611281081103599094690386182078074170366663059556\\
152849646543540637835132530027625323121019017115313612010591786655453628211510631]
\end{gather*}
\vspace{-0.5cm}
\begin{verbatim}
    19 minutes
\end{verbatim}
\begin{verbatim}
  > length(%[5]);
\end{verbatim}
\vspace{-0.3cm}
\begin{equation*}
336
\end{equation*}
How many primes are there in the above defined Smarandache
sequences of triangular numbers? This seems to be an open question.

Another interesting question is to find triangular numbers
in the Smarandache sequences of triangular numbers.
We begin by defining in \textsf{Maple} the boolean function
\texttt{istriangular}.
\begin{verbatim}
  > istriangular := n -> evalb(nops(select(i->evalb(whattype(i)=integer),
  >                                        [solve(t(k)=n)])) > 0):
\end{verbatim}
There exist two triangular numbers (1 and 136) among
the first 1000 terms of the Smarandache consecutive
sequence of triangular numbers;
\begin{verbatim}
  > st := time():
  > select(istriangular,Lscs1000);
  > printf("%a seconds",round(time()-st));
\end{verbatim}
\vspace{-0.3cm}
\begin{equation*}
[1, 136]
\end{equation*}
\vspace{-0.5cm}
\begin{verbatim}
    6 seconds
\end{verbatim}
while the other Smarandache sequences of triangular numbers only show,
among the first 1000 terms, the trivial triangular number 1:
\begin{verbatim}
  > st := time():
  > select(istriangular,Lrss1000);
  > printf("%a seconds",round(time()-st));
\end{verbatim}
\vspace{-0.3cm}
\begin{equation*}
[1]
\end{equation*}
\vspace{-0.5cm}
\begin{verbatim}
    6 seconds
\end{verbatim}
\begin{verbatim}
  > st := time():
  > select(istriangular,Lsms1000);
  > printf("%a seconds",round(time()-st));
\end{verbatim}
\vspace{-0.3cm}
\begin{equation*}
[1]
\end{equation*}
\vspace{-0.5cm}
\begin{verbatim}
    10 seconds
\end{verbatim}
\begin{verbatim}
  > st := time():
  > select(istriangular,Lsss1000);
  > printf("%a seconds",round(time()-st));
\end{verbatim}
\vspace{-0.3cm}
\begin{equation*}
[1]
\end{equation*}
\vspace{-0.5cm}
\begin{verbatim}
    6 seconds
\end{verbatim}

Does exist more triangular numbers in the Smarandache sequences of triangular numbers?
This is, to the best of our knowledge, an open question needing further investigations.
Since checking if a number is triangular is much faster than to check if a number is prime,
we invite the readers to continue our search of triangular numbers
for besides the 1000th term of the Smarandache sequences of triangular numbers.
We look forward to readers discoveries.




\begin{thebibliography}{9}

\bibitem{Gupta} Shyam Sunder Gupta,
Smarandache Sequence of Triangular Numbers,
Smarandache Notions Journal, Vol. 14, 2004, pp.~366--368.

\end{thebibliography}
\end{document}